\documentclass[11pt,a4paper]{article}

\usepackage[T1]{fontenc}
\usepackage[utf8]{inputenc}
\usepackage{mathptmx}       
\usepackage{helvet}         
\usepackage{courier}        
\usepackage{graphicx}        
\usepackage{multicol}        
\usepackage[bottom]{footmisc}
\usepackage{amssymb, amsmath, latexsym}
\usepackage{mathrsfs}
\usepackage{eurosym}
\usepackage{paralist}
\usepackage[margin=1in]{geometry}

\usepackage{amsfonts}
\usepackage{amssymb,enumerate}
\newcommand{\qed}{}
\usepackage{breqn}

\usepackage{relsize}

\allowdisplaybreaks

\makeindex             
\newenvironment{proof}[1][Proof]%
{\vspace{0.5em}\par\textit{{#1}.}\;\;}{\phantom{|}\hfill$\Box$\par}

\begin{document}

\title{Optimal forward contract design for inventory: a~value-of-waiting analysis}

\author{Roy O. Davies\and Adam J. Ostaszewski}
%
%
%
\maketitle

\textbf{Abstract. } A classical inventory problem is studied from the perspective of
embedded options, reducing inventory-management to the design of
optimal contracts for forward delivery of stock (commodity).
Financial option techniques \`{a} la Black-Scholes are invoked to
value the additional `option to expand stock'. A simplified approach
which ignores distant time effects identifies an optimal `time to
deliver' and an optimal `amount to deliver' for a production process
run in continuous time modelled by a Cobb-Douglas revenue function.
Commodity prices, quoted in initial value terms, are assumed to
evolve as a geometric Brownian process with positive drift. Expected
revenue maximization identifies an optimal `strike price' for the
expansion option to be exercised, and uncovers the underlying
martingale in a truncated (censored) commodity price. The paper
establishes comparative statics of the censor in terms of drift and
volatility, and uses asymptotic approximation for a tractable
analysis of the optimal timing.

\section{Problem formulation and model}

We enhance a classical inventory-management problem by studying its
embedded options, reducing the problem to the design of optimal
contracts for forward delivery of inventory. The approach borrows
much from the Black-Scholes model for valuing financial options (see
Musiela and Rutkowski \cite[Chapter 5]{MusR}) and uncovers the underlying
martingale to be a truncated (right-censored) discounted commodity
price.

A production process runs continuously over a unit time interval and
the manager is permitted to acquire raw input materials at two
dates: initially,
at time $t=0,$ and again at one other time $%
\theta <1$, selected freely, but committed to at time $t=0$. This
framework is intended as a proxy for a multi-stage inventory
management problem, since `proximal' effects of forward contracting,
as represented by the date $\theta$, are more significant than any
additional `distal' dates for forward delivery. Distal dates for
additional forward deliveries are thus neglected in this model (see
the `Interpretation' paragraph at the end of Section 4). Inputs
are consumed in a continuous production process which creates an
instantaneous revenue rate at time $t$ equal to $f(x_{t})$ (quoted
in present-value terms), where $x_{t}$ is the instantaneous input
rate of consumed material. To begin with, $f(x)$ is standardly an
\textbf{Inada-type} increasing function, viz. twice differentiable,
unboundedly increasing from zero, with slope unbounded at the origin
and strictly decreasing to zero at infinity; eventually $f(x)$ is
specialized to a Cobb-Douglas production function. The revenue from
any interval $[a,b]$ is
taken to be%
\begin{equation*}
\int_{a}^{b}f(x_{t})dt.
\end{equation*}%
If the manager decides to use up a proportion $\theta x $ in the period $[0,$
$\theta ]$ then, with $\theta $ fixed, the Euler-Lagrange equation implies that
a constant instantaneous input rate equal to $x$ is optimal. A
further quantity $(1-\theta )y$
may similarly be consumed in the remaining time interval. If the quantity $%
(1-\theta )y$ is made up from a contracted forward delivery of
$(1-\theta )u$ and from a possible supplement, purchased at time $\theta$, of a non-negative quantity $%
(1-\theta )z,$ the revenue from the second interval will be%
\begin{equation*}
\int_{\theta }^{1}f(x_{t})dt=(1-\theta )f(u+z).
\end{equation*}

Values here and below are quoted in discounted terms, i.e.
present-value terms relative to time $t=0$. (We side-step a
discussion of the relevant discount factor. In brief, discounting
would be done relative to the required rate of return on capital
given the risk-class of the investment project; see Dixit and
Pindyck \cite[Chapter 4, Section 2]{DixP}.)

Whilst the model of revenue assumes a steady (deterministic) market
for the output,
the input prices are assumed stochastic. (We prefer this modelling choice over the more general approach of including also a stochastic output price. Indeed, what then determines optimal behaviour is the ratio of the two prices; so, in a sense, the present simpler arrangement subsumes it.) Specifically, we suppose that at time $%
0 $ the price of inputs is $b_{0}=1,$ and that, as time $t$
progresses, the present value of the spot price, $b_{t}$, follows
the stochastic differential equation:
\begin{equation}
\frac{db_{t}}{b_{t}}=\bar{\mu}dt+\bar{\sigma}dw_{t},  \label{sde}
\end{equation}%
with $w_{t}$ a standard Wiener process. It is assumed that the
constant growth rate $\bar{\mu}$ is positive, so that the expected
(present-value/discounted) price at time $t$ is $e^{\bar{\mu}t}$;
thus
the price is expected to grow above the initial price of unity. The price $%
b_{t}$ is log-normally distributed with a mean which we denote by $\nu=(\bar{\mu}%
-\frac{1}{2}\bar{\sigma}^{2})t$ and a variance $\sigma
^{2}=\bar{\sigma}^{2}t.$
Write $q_{t}(\cdot)=q(\cdot,\bar{\mu}t,\bar{\sigma}\sqrt{t})$ for the density of $%
b_{t}$. Conditional on the initial choice of $%
\theta $, the expected future revenue consequent on the choice of $x,u$ and $%
z $ (with $z$ selected at time $\theta )$ is%
\begin{equation}
\theta(f(x)-x) +(1-\theta )\left(\int_{0}^{\infty
}\{f(z+u)-bz\}q_{\theta}(b)db-u\right). \label{Total}
\end{equation}

This is a classical inventory problem but amended by the inclusion
explicitly of the `option to expand inventory' (choice of $z$) and
of a `forward' contract (choice of $u$). We will evaluate the
embedded option in a framework reminiscent of Black-Scholes
option-pricing. The `forward contract' is construed here as a
contract signed at the earlier date $t=0$ with an agreed
specified delivered quantity, $u,$ a specified delivery date $%
t=\theta ,$ and a price standardized here to $unity$ per unit
delivered. The latter standardization fixes the unit of money,
since, in the absence of arbitrage and storage costs, as is
well-known, the forward price equals the price of inputs at the
initial time of contracting, compounded up to term-value at the
required rate of interest. Note that the advance purchase of $u$ has
by assumption nil resale value on delivery. This makes the delivered
asset a `non-tradeable' commodity, so that the usual martingale
valuation approach applied to a discounted security price is not
immediately appropriate; our analysis makes recourse to dynamical
programming, as in Eberly and Van Mieghem \cite{EbeM}, and thereby
identifies the underlying martingale structure via an appropriately
truncated (right-censored) price.

Apart from offering a real-options approach with optimal design in
mind, in contrast to the classical inventory literature (see for
instance Bensousssan et al. \cite{Ben}, or Scarf \cite{Sca}), an additional
contribution of the current paper is to provide information about
the sensitivity in regard to model parameters of the critical
`strike price' for stock expansion (its comparative statics and
asymptotics), an issue omitted from consideration in Eberly and Van
Mieghem \cite{EbeM}.

The current study of profit dependence on timing, drift and variance
is motivated by the general discrete-time multi-period model of
Gietzmann and Ostaszewski \cite{GieO2}, but with the simplifying removal
of costly liquidation of inventory. There the latter feature was
necessary for a more comprehensive study into the dependence of a
firm's `future value' on accounting data. Such themes are explored in \cite{Ost} in this volume.

Our option-based analysis is simpler than \cite{ChaCS}, though similar in spirit. There the (retailer's) inventory control problem addresses re-distribution of a storable product; one uses a (long) forward for delivery combined with an option to dispose of any excess (a put, with a lower salvage price) coupled with an option for additional supply (a call, with a penalty cost for `emergency supply'); for background on these `option' terms see e.g. \cite{Hul}. A similar approach, albeit in discrete time, is taken in \cite{KouPQ} using at each date a continuum of puts and calls maturing at the next date taken together with a short (negative) forward.

The rest of the paper is organized as follows. In \S 2 we study optimality conditions, which identify a threshold price level (the price censor) above which it is not worth purchasing the input. We consider its sensitivity (comparative statics) to price drift and volatility in \S 3. Then in \S 4 we assess the expected revenue and in \S 5 the optimal timing. Proofs (sensitivity analysis) is spread across \S 6 and \S 7.

\section{Optimality: the censor and value of waiting}

From (\ref{Total}) the optimization problem separates into maximization of $%
f(x)-x$ (with solution specified by $f^{\prime }(x)=1$) and the
maximization, over choice of scalar $u\geq 0$ and function $z(b)$,
of the (time  $t=0$) expectation
\begin{equation}
E[f(z(b)+u)-bz(b)]-u.  \label{problem}
\end{equation}

\noindent\textbf{Definition.} For any Inada-type, strictly concave function
$f(x)$ define the `indirect profit' (i.e. maximised profit) for a
deterministic
price $b$ by%
\begin{equation}
h(b)=\max_{x>0}[f(x)-bx].  \label{Fenchel}
\end{equation}%
Evidently $h(b)=f(I(b))-bI(b),$ where $I$ traditionally denotes the inverse function to $%
f^{\prime }.$\medskip

\noindent\textbf{Theorem 1 (Optimal forward delivered quantity).}
\textit{In the model setting above, with time }$\theta $ \textit{given, let }%
$\tilde{b}:=\tilde{b}(\mu ,\sigma ,\theta )$\ \textit{be the scalar
solving the equation}%
\begin{equation}
E[b_{\theta }\wedge \tilde{b}]=b_{0}=1,  \label{cen}
\end{equation}%
\textit{where }$b_{\theta }$\textit{\ denotes the random price at time }$%
\theta $\textit{. Then the profit-optimizing level of advance
purchase }$u=u(\mu
,\sigma ,\theta )$\textit{\ for {\rm (\ref{problem})} satisfies }%
\begin{equation}
f^{\prime }(u)=\tilde{b},  \label{hedge}
\end{equation}%
\textit{and the optimal expected profit is given by}%
\begin{equation}
g(\mu ,\sigma )=E[h(b),b\leq \tilde{b}]+h(\tilde{b})\cdot\Pr
[b>\tilde{b}]. \label{profit}
\end{equation}

\begin{proof} With $\beta $ arbitrary, select $u$ with
$\beta =f^{\prime }(u).$ Note that $h(\beta )=f(u)-\beta u$ and
$h^{\prime }(\beta )=-u.$ Define the right-censored random variable
$B_{\theta }=B_{\theta
}(\beta )$ by%
\begin{equation*}
B_{\theta }=b_{\theta }\wedge \beta .
\end{equation*}%
For given price $b$ the quantity $z=z(b)$ which maximizes
$f(u+z)-bz$, is either zero, or satisfies the first-order condition
\begin{equation*}
f^{\prime }(z+u)=b.
\end{equation*}%
In view of the monotonicity of $f'$ we thus have $z(b)=0,$ unless
$b\leq \beta .$ Given that $u$ has been purchased at a price of unity,
the profit, when $b_{\theta }\leq \beta ,$ is $f(u+z)-(b_{\theta
}z(b_{\theta })+u)=h(b_{\theta })+(ub_{\theta }-u).$ Otherwise it is $%
f(u)-u=h(\beta )+u\beta -u.$ Thus the expected profit is%
\begin{equation*}
\Pi (\beta ):=E[h(B_{\theta })+uB_{\theta }-u]=E[h(B_{\theta
})]+uE[B_{\theta }]-u.
\end{equation*}%
Differentiating $\Pi$ with respect to $\beta$, and noting that
\[
dE[h(b_{\theta }\wedge \beta )]/d\beta =h^{\prime }(\beta )\Pr [b_{\theta
}\geq \beta ],
\]
we obtain, after some cancellations in view of
$h^{\prime }(\beta )=-u$, the optimality condition $E[B_{\theta
}]=1$ on $\beta .$ The model assumption that $\bar{\mu}$ is positive
ensures the existence of a solution of equation (\ref{cen}). With
$\beta $ set equal to the solution $\tilde{b}$ of equation
(\ref{cen}) we have $\tilde{b}=f^{\prime }(u),$ i.e. (\ref{hedge}).
\qed
\end{proof}
\bigskip

\noindent\textbf{Definition.} In view of the right-censoring of the price $b$
occurring
under the expectation, we call the solution of (\ref{cen}) the \textbf{censor} $%
\tilde{b}=\tilde{b}(\mu ,\sigma ,\theta )$ at time $\theta .$ This
definition follows Gietzmann and Ostaszewski \cite{GieO1}. The censored
variable is thus a martingale.\medskip

\noindent\textbf{Remark. }It is clear from the proof above that the censor
describes the upper limit of those prices which trigger the exercise
of the option to expand stock. So evidently, $\tilde{b}>1.$
We return in the next section to a consideration of its behaviour.
Whilst this threshold role makes the censor similar to the `optimal
ISD control limit' studied by Eberly and Van Mieghem \cite{EbeM}, their
thresholds correspond to Investing/Staying-put/Disinvesting and are
distinct in respect of the treatment of capital depreciation.\medskip

\noindent\textbf{Proposition 1 (Value of waiting). }
\textit{The expected profit $g(\mu ,\sigma )$ defined in
{\rm (\ref{profit})} obtained by optimal forward contracting is no worse
than
the indirect profit $h(1)$ obtained by only using purchases at initial prices, that is}%
\begin{equation*}
h(1)<g(\mu ,\sigma )=E[h(b_{\theta }\wedge \tilde{b})].
\end{equation*}

\begin{proof} This follows from a simple application of
Jensen's inequality,
as $h(b)$ is strictly convex in $b.$ Indeed, we then have%
\begin{equation*}
h(1)=h(E[b_{\theta }\wedge \tilde{b}])<E[h(b_{\theta }\wedge
\tilde{b})].
\end{equation*}
Of course $-h(b)$ is the Fenchel dual of the strictly concave
function $f$, so $-h(b)$ is strictly concave in $b$ (see \cite[Section 12]{Roc}. In the specific case of $f(x)$ twice
differentiable the asserted convexity follows from $h^{\prime \prime
}(b)=-1/f^{\prime \prime }(I(b)),$ where $I$ denotes, as before, the
inverse function of $f'$.\qed
\end{proof}

\section{Sensitivity: Censor comparative statics}

Assuming an Inada-type production function, for the geometric
Brownian model adopted in respect of price as in (\ref{sde}), the
censor equation (\ref{cen}) which defines $\tilde{b}=\tilde{b}(\mu
,\sigma )$ can be re-written as:
\begin{equation}
1=e^{\mu }\Phi (W-\sigma )+\tilde{b}\Phi (-W). \label{sigw1}
\end{equation}%
Here $\Phi (x)=\int_{-\infty }^{x}\varphi (w)dw$, with $\varphi (w)=%
e^{-\frac{1}{2}w^{2}}/\sqrt{2\pi },\quad$denotes the standard normal
cumulative distribution function, $W=w(\tilde{b})$, and
\begin{equation}
w(b):=\frac{\ln b-\nu}{\sigma },\quad\text{   where }\nu=\mu -\frac{1}{2}%
\sigma ^{2}.  \label{substitution}
\end{equation}

This formulation leads naturally to a further definition.

\smallskip

\noindent\textbf{Definition. }The\textbf{\ normal censor} is the implicit
function $W(\mu ,\sigma )$ defined for $\mu, \sigma >0$ as follows:
\begin{equation}
e^{-\mu }=F(W,\sigma ),\text{    where    }F(W,\sigma ):= \Phi
(W-\sigma )+e^{\sigma W-\frac{1}{2}\sigma ^{2}}\Phi (-W).
\label{sigw}
\end{equation}%
We note that $W$ is well defined since $\partial F/\partial W>0.$ It
is helpful to be aware of the hidden connection between the function
$F$ and the normal hazard rate $H(x)=\varphi(x)/\Phi(-x)$ (or its
reciprocal, the Mills' Ratio) and to use properties of this
function. We refer to Kendall and Stuart (\cite[p.104]{KenS}, or Patel
and Read \cite{PatR} for details. From $\varphi
(\sigma-W )=e^{\sigma W-\frac{1}{2}\sigma ^{2}}\varphi (W)$,
\begin{equation*}
F(W,\sigma )=\varphi (\sigma-W
)\left(\frac{1}{H(\sigma-W)}+\frac{1}{H(W)}\right).
\end{equation*}%
From (\ref{sigw}), $W(\mu ,\sigma )$ is decreasing in $\mu $, as $%
e^{-\mu }$ is decreasing. Less obvious is the fact that $W(\mu
,\sigma )$ is increasing in $\sigma $ since in fact $\partial
W/\partial \sigma >1.$ This
is shown in Section 7, where we deduce the comparative statics of $\widetilde{%
b}(\mu ,\sigma )$ from corresponding properties of $W(\mu ,\sigma
).$ The main results proved there are cited below.\medskip

\noindent \textbf{Theorem 5. }\textit{The censor}\textbf{\
}$\widetilde{b}(\mu ,\sigma )$\textbf{\ }\textit{is decreasing in
the drift and is increasing in the standard deviation.}\medskip

These two properties together suggest the following
result, obtained by setting $\mu=\overline{\mu }\theta$ and
$\sigma=\overline{\sigma }\sqrt{\theta }$, and noting that
(\ref{sigw}) permits arbitrary positive $\theta$.\medskip

\noindent \textbf{Theorem 6. }\textit{The censor }$\bar{b}(\theta )=\widetilde{%
b}(\overline{\mu }\theta ,\overline{\sigma }\sqrt{\theta
})$\textit{\ is either unimodal or increasing on the interval $0 <
\theta <\infty$, according as $\overline{\mu
}\ge\frac{1}{2}\overline{\sigma}^2$ or $\overline{\mu
}<\frac{1}{2}\overline{\sigma}^2$.}

\section{Cobb-Douglas revenue: asymptotic results}

We now assume $f(x)$ is Cobb-Douglas, specifically $f(x)=2\sqrt{x},$
so that the indirect profit defined by (\ref{Fenchel}) is
$h(b)=b^{-1}.$ This choice for the power of $x$ inflicts no loss of
generality, because in the presence of a log-normally distributed
price
any other choice of power is equivalent to a re-scaling of $\bar{\mu},\bar{%
\sigma}.$ Substituting into the definition (\ref{profit}) yields%
\begin{equation}
g(\mu ,\sigma )=e^{(\sigma ^{2}-\mu )}\Phi (W+\sigma )+e^{-\mu -\sigma W+%
\frac{1}{2}\sigma ^{2}}\Phi (-W),  \label{EPF}
\end{equation}%
as the profit per unit time arising after the re-stocking date
$\theta $. We also define the associated function%
\begin{equation*}
\bar{g}(\theta ):=g(\bar{\mu}\theta ,\bar{\sigma}\sqrt{\theta
}),
\end{equation*}%
for $0 < \theta <\infty$ (with some re-sizing of $\bar{\mu}
,\bar{\sigma}$ in mind, as in Proposition 4 of Section 5). To study
these functions we are led to analyze the behaviour of first $W(\mu
,\sigma )$ and then $\bar{W}(t)=_{def}W(\bar{\mu}t,\bar{\sigma}%
\sqrt{t}).$ The following are derived in Section 6.\medskip

\noindent\textbf{Proposition 2.} \textit{ For fixed }$\mu>0$,
\[
W(\mu ,\sigma )=-\frac{\mu }{\sigma }+\frac{1}{2}%
\sigma +o(\sigma ) \quad \textit{as } \sigma \rightarrow 0+.
\]
\medskip

\noindent\textbf{Proposition 3.} \textit{ For fixed }$\mu>0$,
\[
W(\mu ,\sigma )=\sigma -\hat{\mu}-\frac{1}{\sigma -%
\hat{\mu}}\{1+o(1)\} \textit{ as } \sigma \rightarrow \infty,
 \textit{ with }\hat{\mu}=-\Phi^{-1}(e^{-\mu}).
\]
\medskip

From (\ref{EPF}) and standard asymptotic estimates of $\Phi (x)$
(see Abramowitz \& Stegan \cite[Section 7]{AbrS} Theorem
A below is immediate. It also turns out that $\bar{W}(t)$ behaves rather like $%
\pm\sqrt{t}$ (except when $\overline{\mu
}=\frac{1}{2}\overline{\sigma}^2$).\medskip

\noindent\textbf{Theorem A (Asymptotic behaviour of the profit }$g(\mu
,\sigma )$){\bf .}\smallskip

(i) $g=e^{\sigma ^{2}-\mu }+o(1/\sigma )$ \textit{as} $\sigma
\rightarrow \infty ;$\smallskip

(ii) $g=e^{-\mu }+(1-e^{-\mu })\Phi (\mu /\sigma )+o(\sigma )$ \textit{as} $%
\sigma \rightarrow 0+.$

\bigskip

\noindent\textbf{Theorem B (Behaviour of the profit }$\bar{g}(\theta )$
\textbf{at the origin).} \smallskip

\noindent\textit{We have }$\bar{g}^{\prime }(0)=\bar{\sigma}^{2}$\textit{\ so that}%
\begin{equation*}
\bar{g}(\theta )=1+\bar{\sigma}^{2}\theta +o(\theta ).
\end{equation*}

\medskip

\noindent\textbf{Theorem C (Asymptotic behaviour of the profit
}$\bar{g}(\theta )$ \textbf{at infinity).} \medskip

(i) \textit{If }$\bar{\sigma}^{2}<\bar{\mu}$\textit{\ we have as}
$\theta
\rightarrow \infty $\textit{\ that}%
\begin{equation*}
\bar{g}(\theta )=1+o(1/\sqrt{\theta })\rightarrow 1+,
\end{equation*}%
\textit{and }$\bar{g}(\theta )$\textit{\ has a maximum whose
location tends to infinity as }$\bar{\sigma}^{2}\rightarrow\bar{\mu}$\textit{.}\medskip

(ii) \textit{If }$\bar{\mu}\leq $\textit{\ }$\bar{\sigma}^{2}<2\bar{\mu}$%
\textit{\ we have as} $\theta \rightarrow \infty $ \textit{\ that}%
\begin{equation*}
\bar{g}(\theta )=1+e^{(\bar{\sigma}^{2}-\bar{\mu})\theta }+o(1/\sqrt{\theta }%
).
\end{equation*}

(iii) \textit{If }$2\bar{\mu}<\bar{\sigma}^{2}$\textit{\ we have as}
$\theta
\rightarrow \infty $ \textit{\ that}%
\begin{equation*}
\bar{g}(\theta )=e^{(\bar{\sigma}^{2}-\bar{\mu})\theta
}+o(1/\sqrt{\theta }).
\end{equation*}

(iv) \textit{If }$\bar{\sigma}^{2}=2\bar{\mu}$\textit{\ we have as}
$\theta
\rightarrow \infty $ \textit{\ that}%
\begin{equation*}
\bar{g}(\theta )=\frac{1}{4}+e^{\bar{\mu}\theta }\Phi (\sqrt{2\bar{\mu}%
\theta })+o(1/\sqrt{\theta })=\frac{1}{4}+e^{\bar{\mu}\theta }+o(1/\sqrt{%
\theta }).
\end{equation*}\medskip

For the proofs, see \S 6.

Figures 1-4 with a parameter value $\bar{\mu}=0.05$ show the
computed graphs of $\bar{g}$ (bold) alongside the relevant
approximation (faint); Figure 4 shows the first to the right of the
two approximations given in the case (iv).\medskip

\textbf{Interpretation.} Under `myopic management', i.e. in the
absence of forward contracting, for a given re-stocking date $\theta $ the expected profit would be $$
E_{0}[h(b_{\theta })]=E[1/b_{\theta
}]=e^{(\bar{\sigma}^{2}-\bar{\mu})\theta }.$$
The theorem thus
implies that forward contracting advantages lose significance as
variance increases, or as the re-stocking date $\theta $ advances.
This ultimately is our justification for excluding any additional
dates for further forward deliveries.
\begin{figure}[!t] 
\centering\includegraphics[width=0.7\textwidth] {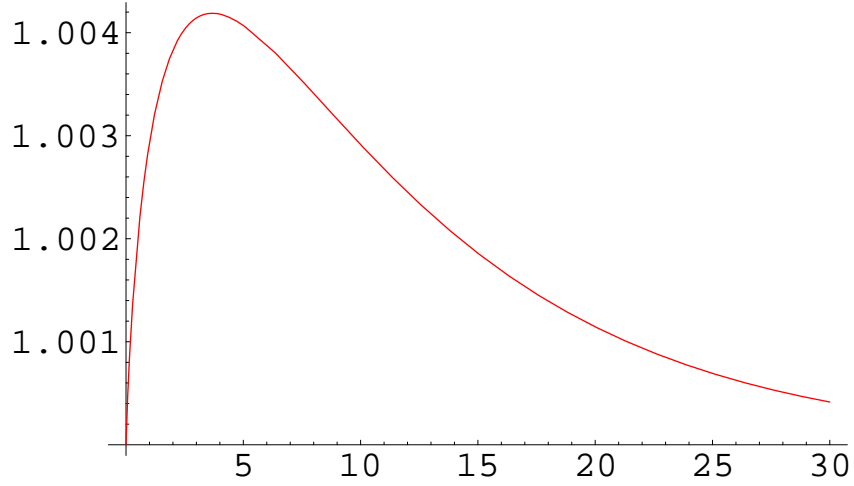}
\caption{Typical graph of $\bar{g}(\theta )$ in the case (i) $\bar{\sigma}%
^{2}<\bar{\mu}$ } \label{fig1}
\end{figure}

\begin{figure}[!t] 
\centering\includegraphics[width=0.7\textwidth] {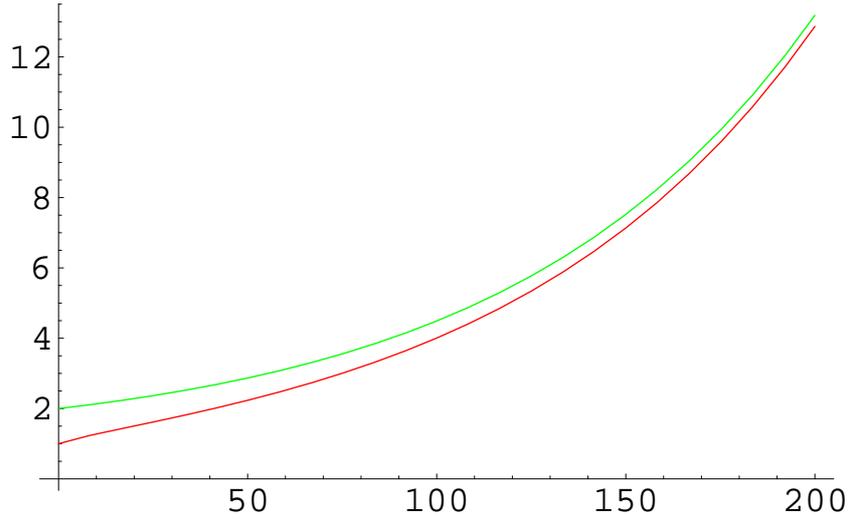}
\caption{Typical graph of $\bar{g}(\theta )$ in the case (ii)
$\bar{\mu}<$ $\bar{\sigma}^{2}<2\bar{\mu}$} \label{fig2}
\end{figure}

\begin{figure}[!t] 
\centering\includegraphics[width=0.7\textwidth] {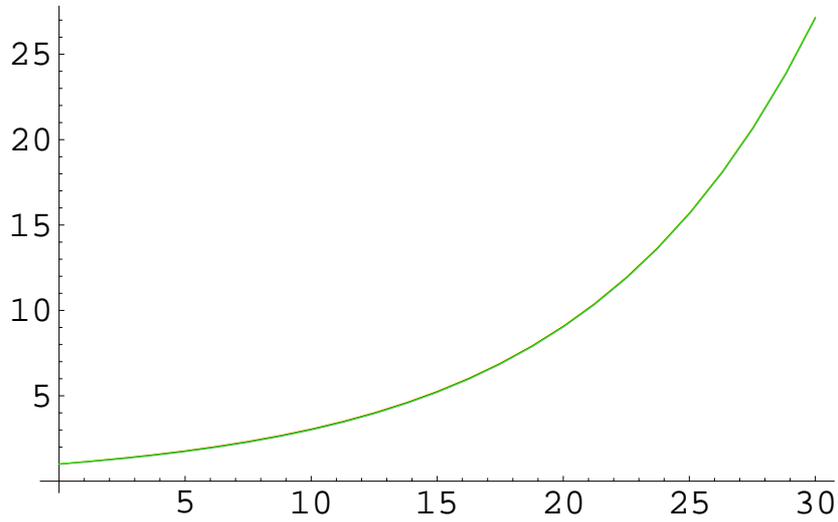}
\caption{Typical graph of $\bar{g}(\theta )$ in the case (iii) $2%
\bar{\mu}<\bar{\sigma}^{2}$} \label{fig3}
\end{figure}

\bigskip

\begin{figure}[!t] 
\centering\includegraphics[width=0.7\textwidth] {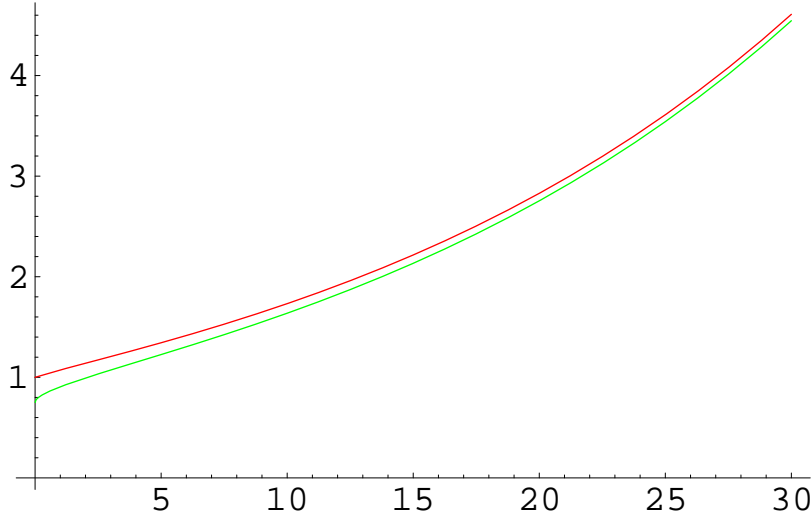}
\caption{Typical graph of $\bar{g}(\theta )$ in the case (iv)
$\bar{\sigma}^{2}=2\bar{\mu}$ and of $\frac{1}{4}+e^{\bar{\mu}\theta }\Phi (\sqrt{2\bar{\mu}%
\theta })$ } \label{fig4}
\end{figure}

\bigskip

\section{Cobb-Douglas optimal timing: estimates}

Assuming as above, again without real loss of generality, that $f(x)=2\sqrt{x},$
we turn now to revenue optimization in respect of the time $\theta $
to be selected
freely in $[0,1].$ Supposing there are no associated management costs in choosing $%
\theta $, the optimal revenue for a selected value of $\theta $ is, from (%
\ref{Total}), given by:%
\begin{equation*}
R(\theta )=\theta +(1-\theta )\bar{g}(\theta ),
\end{equation*}%
the first term being justified by $h(1)=1.$ As $\bar{g}(0)=1$ the
optimal choice of $\theta ,$ assuming such exists, is given by the
following first-order condition:
\begin{equation}
\frac{\bar{g}(\theta )-\bar{g}(0)}{\bar{g}^{\prime }(\theta
)}=1-\theta . \label{FOC}
\end{equation}%

\noindent\textbf{Proposition 4.} \textit{The first-order condition for $R$
in} (\ref{FOC})\textit{ is satisfied for some }$\theta $\textit{\
with} $0<\theta <1.$ \textit{The smallest solution
is a local maximum of }$R$. \textit{If }$\bar{g}$ \textit{is concave on }$%
[0,1],$ \textit{then the solution of }(\ref{FOC})\textit{ is unique.
}

\begin{proof} In general, by Proposition 1 on the Value of Waiting (Section 2), $%
\bar{g}(1)-\bar{g}(0)>0$ and so the first assertion is obvious, as
the right-hand side is zero at $\theta =1$ and is positive at
$\theta =0;$
indeed, by Theorem B above, the left-hand side has a limiting value zero as $%
\theta \rightarrow 0+$ for $\overline{\sigma }>0$. If, however,
$\bar{g}(1)-\bar{g}(0)=0$ (i.e. $h$ failed to be strictly convex),
then since the function $\bar{g}$ is initially increasing for
$\theta
>0,$ $\bar{g}$ has an
internal local maximum at $\bar{\theta}$ for some $\bar{\theta}$ with $0<%
\bar{\theta}<1$ (by the Mean Value Theorem). In this case the
first-order condition for $R$ is satisfied by some
$\theta <\bar{\theta},$ since the left-hand side tends to $+\infty $ as $%
\theta \rightarrow $ $\bar{\theta}$.

Any internal solution $\theta ^{\ast }$ to equation (\ref{FOC}) has
$\bar{g}^{\prime }(\theta ^{\ast })>0$ and so the second assertion
follows since $R^{\prime }(\theta
^{\ast }-)>0$ and $R^{\prime }(\theta ^{\ast }+)<0.$ Observe that if $\bar{g}%
^{\prime \prime }(\theta )<0$, then we have%
\begin{equation*}
\frac{d}{d\theta }\left( \frac{\bar{g}(\theta
)-\bar{g}(0)}{\bar{g}^{\prime }(\theta )}\right) =1-\bar{g}^{\prime
\prime }(\theta )\frac{\bar{g}(\theta )-\bar{g}(0)}{[\bar{g}^{\prime
}(\theta )]^{2}}>0,
\end{equation*}%
so the third assertion is clear; indeed concavity ensures that the
left-hand side of (\ref{FOC}) is an increasing function of $\theta.$
\qed
\end{proof}

One would wish to improve on Proposition 4 to show in more general
circumstances (beyond the concavity which can sometimes fail, as
Figure 1 shows) that (\ref{FOC}) has a unique solution, and to study
dependence on the two parameters of the problem. This appears
analytically intractable. For the purposes of gaining an insight we
propose therefore to replace $\bar{g}(t)$ by a function related to
it through asymptotic analysis (as $t$ varies), on the grounds that
from numeric observation the substitute is qualitatively similar.
Examination of behaviour for large $t$ may be justified by re-sizing
the parameters $\bar{\mu},\bar{\sigma}$ which enables the
termination date to become `large'. This observation then ushers in
the advantages of the asymptotic viewpoint.

Guided by Theorems B and C, we are led to a considerably simpler
problem
obtained by making one of two `typical' substitutions for $\bar{g}(\theta )$, namely
\begin{equation*}
1+A\theta e^{-\alpha \theta},\text{ if
}\bar{\sigma}^{2}<\bar{\mu},\text{ or}\qquad e^{\alpha
\theta},\text{ if }\bar{\mu}<\bar{\sigma}^{2},
\end{equation*}%
according as variance is low, or high. Here $\alpha =|\bar{\sigma}^{2}-\bar{%
\mu}|>0$. The substitution in the first of the two situations fits
qualitatively with numeric observation on the \textit{form} of
$\bar{g}$ (see the Figure 1); it agrees in the second situation with
the general form observed in other Figures and also the asymptotic
form as $t\rightarrow \infty .$\medskip

\textbf{Case (i): $\protect\alpha =\bar{\protect\mu}-\bar{\protect%
\sigma}^{2}>0$.} In this case the optimum time $\theta $ is the
solution of

\begin{equation*}
\theta /(1-\alpha \theta )=1-\theta ,
\end{equation*}%
a quadratic relation, leading to the explicit formula
\begin{equation*}
\theta =\theta (\alpha ):=\frac{1}{2}-\frac{1}{\alpha }\left( -1+\sqrt{1+%
\frac{\alpha ^{2}}{4}}\right) ,
\end{equation*}%
so that, as $\alpha $ increases from zero, the optimal time $\theta
$ recedes from the mid-point towards the origin. That is, low
volatilities bring the replenishment timing back.\medskip

\textbf{Case (ii): $\protect\alpha =\bar{\protect\sigma}^{2}-\bar{%
\protect\mu}>0$.} The first-order condition here reduces to%
\begin{equation*}
(1-e^{-\alpha \theta })/\alpha =1-\theta ,
\end{equation*}%
with a unique solution in the unit interval. Here we can use a
quadratic approximation for the exponential term and solve for
$\theta $ to obtain, for $\alpha<2$, the
approximation%
\begin{equation*}
\theta (\alpha )=\frac{1}{1+\sqrt{1-\alpha /2}},
\end{equation*}%
so that the optimal choice of $\theta $ is close to the midpoint
$\theta =1/2,$ when $\alpha $ is small, but advances, as $\alpha $
increases, towards unity (as a direct computation shows). That is,
high volatilities bring the replenishment position forward (meaning
that waiting longer, beyond the mid-term, is optimal for higher
volatilities).

\section{Asymptotic analysis: the proofs}

In this section we give outline arguments (for further details, see the Appendix) leading to the
Propositions 2 and 3 and Theorems B and C of Section 4.\medskip

\noindent\textbf{Lemma 1. } \textit{We have for fixed} $\mu $%
\begin{equation*}
\lim_{\sigma \rightarrow 0+}W(\mu ,\sigma )=-\infty ,\text{ and }%
\lim_{\sigma \rightarrow 0+}\sigma W(\mu ,\sigma )=-\mu .
\end{equation*}%
This is proved directly from the definition of $W(\mu ,\sigma ).$ We
now prove:\medskip

\noindent\textbf{Proposition 5. } \textit{ For }$\mu>0$,
\[
W(\mu ,\sigma )=-\frac{\mu }{\sigma }+\frac{1}{2}%
\sigma +o(\sigma ) \textit{ as } \sigma \rightarrow 0+.
\]

\begin{proof} For an intuition, note that for small enough
$\sigma $ we have $e^{-\mu }\simeq e^{\sigma W-\frac{1}{2}\sigma
^{2}}$and so
\begin{equation*}
W(\mu ,\sigma )\sim -\frac{\mu }{\sigma }+\frac{1}{2}\sigma .
\end{equation*}%
This argument can be embellished as follows. For any non-zero
$\varepsilon $ let
\begin{equation*}
W(\varepsilon ):=-\frac{\mu }{\sigma }+\frac{1}{2}\sigma +\sigma
\varepsilon ,
\end{equation*}%
so that $\sigma W(\varepsilon )-\frac{1}{2}\sigma ^{2}=-\mu +\sigma
\varepsilon $ and hence
\begin{equation*}
\sigma -W(\varepsilon )=\frac{\mu }{\sigma }+\frac{1}{2}\sigma
-\sigma \varepsilon .
\end{equation*}%
We will prove that for positive $\varepsilon $ we have, for small enough $%
\sigma $, that%
\begin{equation*}
W(-\varepsilon )<W(\mu ,\sigma )<W(\varepsilon ).
\end{equation*}%
This is achieved by showing that for all small enough $\sigma $ the
expression below has the same sign as~$\varepsilon :$%
\begin{equation*}
D(\sigma ):=F(W(\varepsilon ),\sigma )-F(W(\mu ,\sigma ),\sigma
)=F(W(\varepsilon ),\sigma )-e^{-\mu }.
\end{equation*}%
This implies the Proposition. Now $D(0+)=0$ and, since $D(\sigma
)=\Phi (W(\varepsilon )-\sigma )+e^{\sigma W(\varepsilon
)-\frac{1}{2}\sigma ^{2}}\Phi (-W(\varepsilon )),$
\begin{eqnarray*}
D^{\prime }(\sigma ) &=&e^{-\frac{1}{2}(-W(\varepsilon )+\sigma )^{2}}\frac{1%
}{\sqrt{2\pi }}\{-\frac{\mu }{\sigma ^{2}}+\frac{1}{2}-\varepsilon
\}+e^{-\mu +\sigma ^{2}\varepsilon }\{2\sigma \varepsilon \}(1+o(\sigma )) \\
&&+e^{-\mu +\sigma ^{2}\varepsilon }e^{-\frac{1}{2}W(\varepsilon )^{2}}\{%
\frac{\mu }{\sigma ^{2}}+\frac{1}{2}+\varepsilon \}.
\end{eqnarray*}%
Note that the first and third terms contain a factor $\sigma \exp
[-\mu ^{2}/\sigma ^{2}],$ which is small compared to $\sigma .$
So for small
enough $\sigma $ the derivative $D^{\prime }(\sigma )$ has the same sign as $%
\varepsilon .$ So the same is true for $D(\sigma)$. \qed
\end{proof}
\bigskip

\noindent \textbf{Definitions}. Recall from (\ref{sigw}) that $\partial F/\partial W>0$
and $F(-\infty ,\sigma )=0,$ $F(+\infty ,\sigma )=1.$ Let $m$ be
fixed; for
the purposes only of the current section it it convenient to define%
\begin{equation*}
\bar{\Phi}(m)=1-\Phi (m)
\end{equation*}%
and to introduce, also as a temporary measure, a variant form $\hat{W}(m,\sigma )$
of $W(m,\sigma )$ obtained by replacing
$e^{-\mu}$ in (\ref{sigw}) by $\bar{\Phi}(m)$ so that now
\begin{equation}
F(\hat{W}(m,\sigma ),\sigma )=\bar{\Phi}(m)<1.  \label{redefine}
\end{equation}\medskip

\noindent\textbf{Claim. } \textit{For $c$ any constant }
$$
\lim_{\sigma \rightarrow \infty }F(\sigma -c,\sigma )=\bar{\Phi}(c).
$$\smallskip

The proof is routine.\medskip

\textbf{Conclusion from claim.} Notice the consequences for the choices $%
c=(1\pm \varepsilon )m.$ Since
\begin{equation*}
\lim_{\sigma \rightarrow \infty }F(\sigma -(1+\varepsilon )m,\sigma )=\bar{%
\Phi}((1+\varepsilon )m)<\bar{\Phi}(m),
\end{equation*}%
for large enough $\sigma $ we have
\begin{equation*}
F(\sigma -(1+\varepsilon )m,\sigma )<F(W,\sigma ).
\end{equation*}%
Hence for large enough $\sigma $ we have $W>\sigma -(1+\varepsilon
)m.$ Similarly, taking $c=(1-\varepsilon )m$ we obtain $W<\sigma
-(1-\varepsilon )m.$ Thus
\begin{equation*}
W(m,\sigma )=\sigma -m\{1+o(1)\}.
\end{equation*}%
This result can be improved by an argument similar to that of
Proposition 2 by reference to

\begin{equation*}
D(\sigma )=\Phi (\sigma -W)+e^{\sigma W-\frac{1}{2}\sigma ^{2}}\Phi (-W)-%
\bar{\Phi}(m)
\end{equation*}%
to yield the following.

\medskip

\noindent\textbf{Proposition 6. } \textit{With the definition
(\ref{redefine}), for fixed }$m$
\begin{equation*}
\hat{W}(m,\sigma )=\sigma -m-\frac{1}{\sigma -m}\{1+o(1)\}\qquad\text{ (as
}\sigma \rightarrow \infty).
\end{equation*}

\textbf{Conclusion.} $\hat{W}(m,\sigma )=W(\mu ,\sigma )$ when $m=\hat{\mu}$ where $e^{-\mu }=\bar{\Phi}(m).$
Restating this equation as $e^{-\mu }=1-\Phi (\hat{\mu})=\Phi
(-\hat{\mu}),$
we see that $\hat{\mu}>0$ if and only if $\mu >\ln 2,$ since $\hat{\mu}%
=-\Phi ^{-1}(e^{-\mu });$ in particular for small $\mu $ we thus have $\hat{%
\mu}<0.$\medskip

\noindent\textbf{Lemma 2. }
\begin{equation*}
\lim_{\theta \rightarrow 0+}\sqrt{\theta }\bar{W}(\theta )=0\text{ and }%
\lim_{\theta \rightarrow 0+}\bar{W}(\theta )=+\infty \text{ \textit{ for fixed%
} }\bar{\mu},\bar{\sigma}>0.
\end{equation*}\smallskip

This follows again by a routine argument starting from (\ref{sigw}),
but requires the claim below.\medskip

\noindent\textbf{Claim. }
\begin{equation*}
L=\lim_{\theta \rightarrow 0+}\sigma V(\theta )=0.
\end{equation*}\smallskip

The proof here is by contradiction from (\ref{sigw}), assuming $L$
non-zero.\medskip

\noindent\textbf{Proof of Theorem B}. Differentiation of (\ref{sigw}) with respect to $\theta$ gives%
\begin{align*}
-\bar{\mu}e^{-\bar{\mu}\theta }=&\;\varphi (W-\sigma )(W^{\prime
}-\sigma ^{\prime })+e^{\sigma W-\frac{1}{2}\sigma ^{2}}\varphi
(-W)(-W^{\prime
})\\
&\;+\Phi (-W)e^{\sigma W-\frac{1}{2}\sigma ^{2}}[-\frac{1}{2}\bar{\sigma}%
^{2}+(\sigma W)^{\prime }].
\end{align*}%
Now
\begin{align*}
\varphi (W-\sigma )\sigma ^{\prime }=&\;e^{\sigma W-\frac{1}{2}\sigma
^{2}}\varphi (W)\frac{\bar{\sigma}}{2\sqrt{\theta }}\\
=&\;\left( e^{\sigma W-%
\frac{1}{2}\sigma ^{2}}\frac{\varphi (W)}{W}\right) \frac{1}{\theta }\frac{W%
\bar{\sigma}\sqrt{\theta }}{2}\rightarrow \bar{\mu}\cdot 0=0\qquad(\theta \rightarrow 0+),
\end{align*}%
using (\ref{sigw}) and $\lim_{W\rightarrow +\infty
}\varphi
(W)/(W\Phi (-W))=1$ to deal with the bracketed term. Thus%
\begin{equation*}
-\bar{\mu}=\lim_{\theta \rightarrow 0+}[\Phi (-W)(\sigma W)^{\prime
}].
\end{equation*}%
Differentiation of (\ref{EPF}) with respect to $\theta$ gives%
\begin{eqnarray*}
\bar{g}^{\prime } &=&[\bar{\sigma}^{2}-\bar{\mu}]e^{(\sigma ^{2}-\mu
)}\Phi (W+\sigma )+e^{(\sigma ^{2}-\mu )}\varphi (W+\sigma
)(W^{\prime }+\sigma
^{\prime }) \\
&&+e^{-\mu -\sigma W+\frac{1}{2}\sigma ^{2}}\varphi (-W)(-W^{\prime
})+e^{-\mu -\sigma W+\frac{1}{2}\sigma ^{2}}\Phi (-W)[\frac{1}{2}\bar{\sigma}%
^{2}-\bar{\mu}-(\sigma W)^{\prime }].
\end{eqnarray*}%
Now%
\begin{equation*}
\bar{g}^{\prime }=[\bar{\sigma}^{2}-\bar{\mu}]-\lim_{\theta
\rightarrow 0+}\Phi (-W)[(\sigma W)^{\prime }]=\bar{\sigma}^{2}.
\end{equation*}{\hfill$\Box$\par}

\noindent\textbf{Lemma 3. }If $\frac{1}{2}\bar{\sigma}^{2}\neq $ $\bar{\mu}$,
then \textbf{\ }
\begin{equation*}
\lim_{\theta \rightarrow \infty }\bar{W}(\theta )=\pm \infty .
\end{equation*}\smallskip

\noindent\textbf{Remark. }This leaves the identification of the appropriate
sign as a separate task. The proof is by contradiction from
(\ref{sigw}) by reference to simple properties of the normal hazard
rate $H(w)=\phi (w)/\Phi (-w)$. Lemma 4 below is proved by
contradiction. Lemma 5 clarifies the cross-over case.\medskip

\noindent\textbf{Lemma 4.} $\lim_{\theta \rightarrow \infty }\bar{W}(\theta
)-\sigma =-\infty .$\medskip

\noindent\textbf{Lemma 5. }\textit{If }$\frac{1}{2}\bar{\sigma}^{2}=\bar{\mu},$%
\textit{\ then} $\lim_{\theta \rightarrow \infty }\bar{\sigma}\sqrt{\theta }%
\bar{W}(\theta )=\log 2.$\medskip

\noindent\textbf{Conclusion 1. }\textit{If}\textbf{\ }$\lim_{\theta
\rightarrow
\infty }\bar{W}(\theta )=-\infty ,$ \textit{then for }$\bar{\mu}>\frac{1}{2}$%
\begin{equation*}
\bar{W}(\theta )=-\frac{\bar{\mu}-\frac{1}{2}\bar{\sigma}^{2}}{\bar{\sigma}}%
\sqrt{\theta }+o(\sqrt{\theta }).
\end{equation*}

\noindent\textbf{Lemma 6. }\textit{If }$\bar{\sigma}^{2}<2\bar{\mu},$\textit{\ then}%
\begin{equation*}
\lim_{\theta \rightarrow \infty }e^{-\sigma \bar{W}+(\frac{1}{2}\bar{\sigma}%
^{2}-\bar{\mu})\theta }=1.
\end{equation*}\smallskip

This follows directly from (\ref{sigw}) and Lemmas 3 and 4.
\medskip

\noindent\textbf{Proof of Theorem C}. Lemma 6 establishes case (ii) of Theorem C. Next we note:
\smallskip

\noindent\textbf{Conclusion 2. }\textit{If}\textbf{\ }$\lim_{\theta
\rightarrow \infty }\bar{W}(\theta )=+\infty ,$ \textit{then for }$\bar{\mu}<\frac{1}{2}\bar{\sigma}^{2}$
\begin{equation*}
\bar{W}(\theta )=(\bar{\sigma}-\sqrt{2\bar{\mu}})\sqrt{\theta }+O(1/\sqrt{%
\theta }).
\end{equation*}

Case (iii) of Theorem C follows from this estimate. Combining (ii) and (iii) gives (i). Turning to case (iv),  if $\bar{\sigma}^{2}=2%
\bar{\mu},$ then as $\theta \rightarrow \infty $ we have
$\sigma +\bar{W}(\theta )\rightarrow +\infty $ , by Lemma 5, so since
$\lim_{\theta \rightarrow \infty }e^{\sigma \bar{W}(\theta )}=2$,
and appealing to standard
asymptotic estimates of $\Phi (x)\sim 1-\varphi(x)/x$,  for large $x$,%
\begin{eqnarray*}
\bar{g}(\theta ) &=&e^{(\bar{\sigma}^{2}-\bar{\mu})\theta }\Phi (\sigma +%
\bar{W}(\theta ))+e^{-\sigma \bar{W}+(\frac{1}{2}\bar{\sigma}^{2}-\bar{\mu}%
)\theta }\Phi (-\bar{W}(\theta )) \\
&=&e^{\bar{\mu}\theta }\Phi (\sigma +\bar{W}(\theta ))+e^{-\sigma \bar{W}%
}\Phi (-\bar{W}(\theta )) \\
&=&e^{\bar{\mu}\theta }+\frac{1}{4}+o(1/\sqrt{\theta }).
\end{eqnarray*}{\hfill$\Box$\par}

\section{Censor comparative statics: reprieve}

This section considers the sensitivity of $\widetilde{b}(\mu ,\sigma
)$ to $\mu$ and $\sigma$, and the dependence of $\bar{b}(\theta )=\widetilde{%
b}(\overline{\mu }\theta ,\overline{\sigma }\sqrt{\theta })$ on
$\theta$ as given in Section 3.
\medskip

\noindent\textbf{Theorem 2. }\textit{The censor }$\widetilde{b}(\mu ,\sigma )$\textit{%
\ is decreasing in the drift }$\mu$.

\begin{proof} The derivative of $\tilde{b}%
=\exp (\sigma W+\mu -\frac{1}{2}\sigma ^{2})$ with respect to
$\mu$ is positive iff:%
\begin{equation}
-\sigma \frac{\partial W(\mu ,\sigma )}{\partial \mu }>1.
\label{wmu}
\end{equation}
But differentiation of (\ref{sigw}) and
$$\varphi (W(\mu ,\sigma )-\sigma )=e^{\sigma W-\frac{1}{2}\sigma
^{2}}\varphi (W(\mu ,\sigma ))$$ yield
\begin{equation*}
1=\widetilde{b}\Phi (-W(\mu ,\sigma ))\left( -\sigma \frac{\partial W}{%
\partial \mu }\right) .
\end{equation*}%
So (\ref{wmu}) holds iff $\widetilde{b}\Phi (-W(\mu ,\sigma ))<1.$
But the latter follows from (\ref{sigw1}).
\qed
\end{proof}
\bigskip

\noindent \textbf{Theorem 3. }\textit{The censor }$\widetilde{b}(\mu
,\sigma )$\textit{\ is increasing in the standard deviation}
$\sigma$.

\begin{proof} Differentiating $\tilde{b}%
=\exp (\sigma W(\mu ,\sigma )+\mu -\frac{1}{2}\sigma ^{2})$ with
respect to $\sigma$ yields
\begin{equation*}
\frac{\partial \widetilde{b}}{\partial \sigma }=\widetilde{b}(\mu
,\sigma )\left\{ \sigma \frac{\partial W}{\partial \sigma }+W(\mu
,\sigma )-\sigma \right\} .
\end{equation*}%
Differentiating also the censor equation (\ref{sigw}) with respect to $%
\sigma $, we obtain after some cancellations
\begin{equation*}
\varphi (W(\mu ,\sigma )-\sigma )=e^{\sigma W(\mu ,\sigma )-\frac{1}{2}%
\sigma ^{2}}\Phi (-W(\mu ,\sigma ))\left\{ W(\mu ,\sigma )+\sigma \frac{%
\partial W}{\partial \sigma }-\sigma \right\} .
\end{equation*}%
The bracketed term appearing here and earlier is thus positive, and so $$\partial \widetilde{b}%
(\mu ,\sigma )/\partial \sigma >0.$$\qed
\end{proof}

Using $\varphi
(\sigma-W )=e^{\sigma W-\frac{1}{2}\sigma ^{2}}\varphi (W)$ (cf. Section 2) we note
the identity
\begin{equation}
W(\mu ,\sigma )+\sigma \frac{\partial W}{\partial \sigma }-\sigma =\frac{%
\varphi (W(\mu ,\sigma ))}{\Phi (-W(\mu ,\sigma ))}=H(W(\mu ,\sigma
)), \label{haz}
\end{equation}%
where $H(x)$ denotes the normal hazard rate ($\varphi (x)/\Phi
(-x)$). Since $H(x)>x$ for all $x,$ equation (\ref{haz}) gives
$\partial W/\partial \sigma >1$ for $\sigma>0$. Recalling from
Section 2 that $\partial W/\partial \mu<0$, we have the two
results:

\bigskip

\noindent\textbf{Theorem 4. }\textit{The two functions \ \ }$\sigma W(\mu ,\sigma )-%
\frac{1}{2}\sigma ^{2},\quad W(\mu ,\sigma )-\sigma $\textit{\ are
increasing in $\sigma$ for $\sigma>0$.}

\bigskip

\noindent\textbf{Theorem 5. }\textit{The normal censor }$W(\mu ,\sigma
)$\textit{\ is increasing in standard deviation and decreasing with
drift.}

\medskip
Our final result is the following.
\medskip

\noindent\textbf{Theorem 6. }\textit{The censor }$\bar{b}(\theta )=\widetilde{%
b}(\overline{\mu }\theta ,\overline{\sigma }\sqrt{\theta
})$\textit{\ is either unimodal or increasing on the interval $0 <
\theta < \infty$, according as $\overline{\mu
}\ge\frac{1}{2}\overline{\sigma}^2$, or as $\overline{\mu
}<\frac{1}{2}\overline{\sigma}^2$.}

\begin{proof} Using $\bar{b}\varphi (W)=e^{\mu }\varphi
(W-\sigma )$ and
applying the Chain Rule to $\bar{b}(\theta )=\widetilde{b}(\overline{\mu }%
\theta ,\overline{\sigma }\sqrt{\theta }),$ we obtain%
\[
\theta \Phi (-W)\frac{d\bar{b}(\theta )}{d\theta }=-\mu \{e^{\mu
}\Phi (W-\sigma )\}+\frac{1}{2}\sigma \bar{b}\varphi (W).
\]%
The stationarity condition for $\bar{b}(\theta )$ can be written
using the normal hazard rate $H(x)=\varphi(x)/\Phi(-x)$ as
\begin{equation}
\mu =\frac{1}{2}\sigma H(-W(\mu ,\sigma )+\sigma ),
\label{unimodal}
\end{equation}%
where $\mu =\overline{\mu }\theta $ and $\sigma =\overline{\sigma }\sqrt{%
\theta }$, and $W(\mu ,\sigma )$ is the normal censor as in (\ref{sigw}%
).

We now regard $\mu $ and $\sigma $ as free variables and let $\kappa :=%
\overline{\mu }/\overline{\sigma }^{2}$ be the dispersion parameter.
In this setting we seek a stationary point $\theta $
of $\bar{b}(\theta )$ by first finding the values $\mu=\mu^{\ast}$
and $\sigma=\sigma^{\ast}$ which satisfy the equation
(\ref{unimodal}) simultaneously with
the equation:%
\begin{equation}
\mu =\kappa \sigma ^{2}.  \label{parabola}
\end{equation}%
We will show that this is possible (uniquely) if and only if $\kappa
\geq 1/2$ (i.e. $\overline{\mu}\ge\frac{1}{2}\overline{\sigma}^2$).
Thus for $\overline{\sigma }^{2}>2\overline{\mu }$ the
function $\bar{b}(\theta )$ is increasing, but otherwise has a unique maximum at $\theta =\mu^{\ast} /\overline{%
\mu }={\sigma^{\ast}}^2/\overline{\sigma}^2$.

We begin by noting that (\ref{unimodal}) defines an implicit
function $\mu =\mu (\sigma )$ for all $\sigma >0.$ Indeed,
elimination of $\mu$ between (\ref{sigw}) and (\ref{unimodal}) leads
to
\begin{equation}
\exp \left(-\frac{1}{2}\sigma H(-w+\sigma )\right)=F(w,\sigma ),
\label{omega}
\end{equation}%
and then routine analysis shows that there is a unique solution
$w=\omega (\sigma )$ of (\ref{omega}).
Since $\partial W/\partial \mu <0$, we may recover $\mu (\sigma )>0$, for $%
\sigma >0$, from $\omega (\sigma ):=W(\mu (\sigma ),\sigma )$.

Linearization of both sides of (\ref{omega}) around $\sigma =0,$
yields the equation
\begin{equation*}
H(-w)=2(\varphi(w)+w \Phi(-w)),
\end{equation*}%
with unique solution $w=\omega (0)=0.$ Hence $\lim_{\sigma
\rightarrow 0}W(\mu (\sigma ),\sigma )=0$ and so, for small
$\sigma$, we have the approximation to (\ref{unimodal}) given by the
convex function
\[
\mu (\sigma):=\frac{1}{2}\sigma H(\sigma ).
\]%
Numerical investigation of the positive function $w=\omega (\sigma
)$ finds
its maximum to be 0.051 for $\sigma $ approximately 2.547. To see why, rewrite (%
\ref{omega}) in the equivalent form:%
\[
\exp \left(\frac{1}{2}\sigma ^{2}-\frac{1}{2}\sigma H(-w+\sigma
)\right)=\frac{\Phi (w-\sigma )}{\varphi (\sigma )\sqrt{2\pi
}}+e^{\sigma w}\Phi (-w).
\]%
For fixed $0\leq w\leq 1$, and large $\sigma $, the left-hand side
is close to $e^{\frac{1}{2}(\sigma w-1)}$, in view of the asymptotic
over-approximation $\left(x+1/x\right)$ for $H(x)$ (when
$x$ is large), whereas the first term on the right is of the order of $%
1/(\sigma \sqrt{2\pi })$. Neglecting the latter, and replacing $\Phi
(-w)$ by $\frac{1}{2}$, the solution for $w$ may be estimated by
$(2\log 2-1)/\sigma $.

Finally, using the same asymptotic approximation for $H(\sigma)$, we
may over-approximate $\frac{1}{2}\sigma H(\sigma -\omega (\sigma
))$ by $\frac{1}{2}\sigma ^{2}+\frac{1}{2}$. From here we may conclude that, for $%
\kappa >\frac{1}{2}$, the equations (\ref{unimodal}) and
(\ref{parabola})
have a solution with a crude over-estimate for $\sigma^{\ast}$ given by:%
\[
\sigma ^{2}=\frac{1}{2\kappa -1}.
\]

The supporting line $\mu=\frac{1}{2}H(0)\sigma$ provides the crude
under-estimate $\sigma=1/\kappa\sqrt{2\pi}$. For the special case
$\kappa =\frac{1}{2}$ the solution to (\ref{unimodal}) and
(\ref{parabola}) is $\sigma^{\ast} =4.331$.  For $\kappa
<\frac{1}{2}$ there is no solution, since $\frac{1}{2}\sigma
H(\sigma )>\kappa \sigma ^{2}$ for $\sigma>0.$\qed
\end{proof}
\bigskip
\noindent\textbf{Acknowledgement}. It is a pleasure to thank Alain Bensoussan for his very helpful advice and encouragement.
\bigskip

\noindent\textbf{Postscript}. Harold Wilson (1916-1995, prime minister 1964-70 and 1974-76) famously always emphasized the importance of keeping his options open.

\section*{Appendix}
\noindent \textbf{Proof of Proposition 6}. For convenience put
\begin{equation*}
R(W,\sigma ):=\sqrt{2\pi }F(W,\sigma )=\int_{-W+\sigma }^{\infty }e^{-\frac{1%
}{2}x^{2}}dx+e^{\sigma W-\frac{1}{2}\sigma ^{2}}\int_{W}^{\infty }e^{-\frac{1%
}{2}x^{2}}dx.
\end{equation*}%
Consider an arbitrary non-zero $\varepsilon ;$ let $W_{\varepsilon }:=\sigma
-m-\delta $ and put
\begin{equation*}
\delta :=\frac{1-\varepsilon }{\sigma -m}.
\end{equation*}%
Now, with $D$ as in \S 6,
\begin{eqnarray*}
D(\sigma ) &=&\left( \int_{-W+\sigma }^{\infty }e^{-\frac{1}{2}%
x^{2}}dx+e^{\sigma W-\frac{1}{2}\sigma ^{2}}\int_{W}^{\infty }e^{-\frac{1}{2}%
x^{2}}dx\right) -\int_{m}^{\infty }e^{-\frac{1}{2}x^{2}}dx \\
&=&\left( \int_{m+\delta }^{\infty }e^{-\frac{1}{2}x^{2}}dx-\int_{m}^{\infty
}e^{-\frac{1}{2}x^{2}}dx\right) +e^{\sigma (\sigma -m-\delta )-\frac{1}{2}%
\sigma ^{2}}\int_{\sigma -m-\delta }^{\infty }e^{-\frac{1}{2}x^{2}}dx \\
&=&-\delta e^{-\frac{1}{2}(m+\delta )^{2}}+O(\delta ^{2})+e^{\sigma (\sigma
-m-\delta )-\frac{1}{2}\sigma ^{2}}\frac{1}{\sigma -m-\delta }e^{-\frac{1}{2}%
(m+\delta -\sigma )^{2}}\{1+O(\frac{1}{\sigma ^{2}})\} \\
&=&-\delta e^{-\frac{1}{2}(m+\delta )^{2}}+O(\delta ^{2})+\frac{1}{\sigma
-m-\delta }e^{-\frac{1}{2}(m+\delta )^{2}}\{1+O(\frac{1}{\sigma ^{2}})\} \\
&=&\left( \frac{1}{\sigma -m-\delta }-\delta \right) e^{-\frac{1}{2}%
(m+\delta )^{2}}+O(\delta ^{2})+O(\frac{1}{\sigma ^{2}}) \\
&=&\left( \frac{1}{(\sigma -m)-\frac{1-\varepsilon }{\sigma -m}}-\frac{%
1-\varepsilon }{\sigma -m}\right) e^{-\frac{1}{2}(m+\delta )^{2}}+O(\delta
^{2})+O(\frac{1}{\sigma ^{2}}) \\
&=&\left( \frac{(\sigma -m)^{2}-(1-\varepsilon )\{(\sigma
-m)^{2}-(1-\varepsilon )\}}{(\sigma -m)^{3}-(1-\varepsilon )(\sigma -m)}%
\right) e^{-\frac{1}{2}(m+\delta )^{2}}+O(\frac{1}{\sigma ^{2}}) \\
&=&\frac{\varepsilon (\sigma -m)^{2}+(1-\varepsilon )^{2}}{(\sigma
-m)^{3}-(1-\varepsilon )(\sigma -m)}e^{-\frac{1}{2}(m+\delta )^{2}}+O(\frac{1%
}{\sigma ^{2}}) \\
&=&\frac{\varepsilon }{\sigma -m}e^{-\frac{1}{2}(m+\delta )^{2}}+O(\frac{1}{%
\sigma ^{2}}),
\end{eqnarray*}%
and this has the same sign as $\varepsilon .$ Thus, for $\varepsilon >0$,
\begin{equation*}
R(W_{-\varepsilon },\sigma )<R(W(m,\sigma ),\sigma )<R(W_{\varepsilon
},\sigma ),
\end{equation*}%
and so, since $\partial R(W,\sigma )/\partial W>0,$
\begin{equation*}
W_{-\varepsilon }<W(m,\sigma )<W_{\varepsilon }.\qquad \square
\end{equation*}
\bigskip
\noindent \textbf{Proof of Lemma 2. }We begin with the associated Claim (%
\S 6 above), for which we need first to put
\begin{equation*}
V:=V(\theta )=\bar{W}(\theta )-\bar{\sigma}\sqrt{\theta },
\end{equation*}%
and then to note (by the definition of the normal sensor in \S 3):
\begin{equation*}
(e^{-\mu }-1)-\{\Phi (-\sigma -V)-\Phi (-V)\}=[e^{\sigma V+\frac{1}{2}\sigma
^{2}}-1]\Phi (-\sigma -V).\tag{*}
\end{equation*}%
From here, for some $V^{\ast }$ between $V$ and $V+\sigma ,$%
\begin{equation*}
(e^{-\mu }-1)-\sigma \varphi (V^{\ast })=[e^{\sigma V+\frac{1}{2}\sigma
^{2}}-1]\Phi (-\sigma -V),
\end{equation*}
so that
\begin{equation*}
-\bar{\mu}\theta +\sigma \varphi (V^{\ast })\sim \lbrack e^{\sigma V+\frac{1%
}{2}\sigma ^{2}}-1]\Phi (-\sigma -V).
\end{equation*}

\noindent \textit{Proof of Claim. }Suppose $L=\lim_{\theta \rightarrow
0+}\sigma V(\theta )\neq 0$ along a sequence of values of $\theta ;$ then
\begin{equation*}
V(\theta )\approx L/(\bar{\sigma}\sqrt{\theta }):\qquad \sigma \varphi
(V^{\ast })\sim \bar{\sigma}\sqrt{\theta }\exp (-L^{2}/\bar{\sigma}%
^{2}\theta )/\sqrt{2\pi }
\end{equation*}%
and so
\begin{equation*}
-\bar{\mu}\theta \{1-(\bar{\sigma}/\bar{\mu}\sqrt{\theta })\exp (-L^{2}/\bar{%
\sigma}^{2}\theta )/\sqrt{2\pi }\sim -\bar{\mu}\theta .
\end{equation*}%
So, for small enough $\theta ,$
\begin{equation*}
\lbrack e^{\sigma V+\frac{1}{2}\sigma ^{2}}-1]\Phi (-\sigma -V)<0,
\end{equation*}%
so that $\bar{V}\leq 0.$ Suppose first that $\bar{V}=-\infty ;$ then $L=0$,
since $\Phi (\infty )=1$ reduces equation (*) to
\begin{equation*}
0=(e^{L}-1),
\end{equation*}%
contradicting $L\neq 0$.
Likewise, from equation (*), the finiteness of $\bar{V}$ yields
$L=0$, a final contradiction. $\square _{\text{claim}}$

Turning to Lemma 2 proper, put $\bar{V}:=\lim_{\theta \rightarrow 0+}V(\theta ).
$ As above
\begin{equation*}
(e^{-\mu }-1)+\sigma \varphi (V^{\ast })\sim \lbrack e^{\sigma V+\frac{1}{2}%
\sigma ^{2}}-1]\Phi (-\sigma -V).
\end{equation*}

By the Claim, $\sigma V$ is small; so we may expand the exponential and,
dividing by $\sigma =\bar{\sigma}\sqrt{\theta },$ obtain
\begin{equation*}
-\frac{\bar{\mu}}{\bar{\sigma}}\sqrt{\theta }+\varphi (V^{\ast })=(V+\frac{1%
}{2}\sigma )\Phi (-\sigma -V).
\end{equation*}%
If $V\rightarrow \bar{V}$ a finite limit, then the Mills ratio (hazard
rate), defined by
\begin{equation*}
H(\bar{V}):=\frac{\varphi (\bar{V})}{\Phi (-\bar{V})},
\end{equation*}%
satisfies $H(\bar{V})=\bar{V},$ a contradiction, since the ratio is always
greater than $\bar{V}.$ Thus the limit $\bar{V}$ must be infinite, hence $%
\varphi (\bar{V})=0.$ So $\bar{V}=+\infty ,$ as otherwise $\bar{V}=-\infty $
leads to the contradiction
\begin{equation*}
0=\bar{V}\Phi \left( -\bar{V}\right) =\bar{V}\cdot 1.\qquad \square
\end{equation*}%
\noindent \textbf{Proof of Lemma 3. }As in the definition of the normal
censor
\begin{equation*}
e^{-\bar{\mu}\theta }=\Phi (\bar{W}(\theta )-\sigma )+e^{\sigma \bar{W}%
(\theta )-\frac{1}{2}\sigma ^{2}}\Phi (-\bar{W}(\theta )),
\end{equation*}%
or%
\begin{equation}
e^{-\bar{\mu}\theta -\sigma \bar{W}(\theta )+\frac{1}{2}\sigma
^{2}}=e^{-\sigma \bar{W}(\theta )+\frac{1}{2}\sigma ^{2}}\Phi (\bar{W}%
(\theta )-\sigma )+\Phi (-\bar{W}(\theta )),  \tag{**}
\end{equation}%
\begin{equation*}
e^{-\bar{\mu}\theta -\sigma \bar{W}(\theta )+\frac{1}{2}\sigma ^{2}}=\Phi (-%
\bar{W}(\theta ))+\varphi (\bar{W}(\theta ))/H(\sigma -\bar{W}(\theta )),
\end{equation*}%
where, as above, $H(.)$ denotes the hazard rate. Assume that $\bar{W}(\theta
)\rightarrow \bar{w}.$ We are to prove that $\bar{w}$ is not finite. We
argue by cases.

\noindent Case 1: $\frac{1}{2}\bar{\sigma}^{2}>\bar{\mu}.$ The left hand side is unbounded,
whereas the right-hand side is bounded for large~$\theta $ by
\begin{equation*}
1+\varphi (\bar{w})/(\bar{\sigma}\sqrt{\theta }-\bar{w}).
\end{equation*}%
\noindent Case 2: $\frac{1}{2}\bar{\sigma}^{2}<\bar{\mu}$. Letting $\theta \rightarrow
\infty $ gives the contradiction:
\begin{equation*}
0=\Phi (-\bar{w})+0.\qquad \square
\end{equation*}
\bigskip
\noindent \textbf{Proof of Lemma 4.} As before, if $V:=\bar{W}(\theta
)-\sigma ,$ then%
\begin{equation*}
e^{-\mu }=\Phi (V)+e^{\sigma V+\frac{1}{2}\sigma ^{2}}\Phi (-\sigma -V).
\end{equation*}%
Suppose $V\rightarrow -\infty $ is false. Then either $V\rightarrow \infty ,$
or $V\rightarrow \bar{V},$ a finite limit. In either case we have
\begin{equation*}
e^{\sigma V+\frac{1}{2}\sigma ^{2}}\Phi (-\sigma -V)\leq e^{-\frac{1}{2}%
V^{2}}\varphi (V+\sigma )/(V+\sigma )\rightarrow 0,
\end{equation*}%
as $\theta \rightarrow \infty $ (since $\sigma \rightarrow \infty ).$ This
implies that $0=\Phi (\bar{V}),$ a contradiction in either case. So $%
V\rightarrow -\infty .$ $\square $
\bigskip

\noindent \textbf{Proof of Lemma 5. }As before suppose $\bar{W}(\theta
)\rightarrow \bar{w}.$ If $\bar{w}<0$ (possibly infinite), then we have in
the limit $\Phi (-\bar{w})=\infty ,$ a contradiction. If $0<\bar{w}<\infty ,$
then by (**) above $0=\Phi (-\bar{w}),$ again a contradiction. This leaves two
possibilities: either $\bar{w}=\infty $ or $\bar{w}=0.$

Suppose the former. Noting that%
\begin{equation*}
1=\lim_{\theta \rightarrow \infty }[e^{\bar{\mu}\theta }\Phi (\bar{W}(\theta
)-\sigma )+e^{\sigma \bar{W}(\theta )}\Phi (-\bar{W}(\theta ))],
\end{equation*}%
then $e^{\sigma \bar{W}(\theta )}\Phi (-\bar{W}(\theta ))$ is bounded. But%
\begin{equation*}
\lim_{\theta \rightarrow \infty }e^{\sigma \bar{W}(\theta )}\Phi (-\bar{W}%
(\theta ))=\lim_{\theta \rightarrow \infty }\frac{e^{\sigma \bar{W}(\theta
)}e^{-\frac{1}{2}\bar{W}^{2}}}{\bar{W}(\theta )\sqrt{2\pi }}=\lim_{\theta
\rightarrow \infty }e^{[\bar{W}(\theta )(\sigma -W)]}\frac{e^{+\frac{1}{2}%
\bar{W}^{2}}}{\bar{W}(\theta )\sqrt{2\pi }}=\infty ,
\end{equation*}%
by Lemma 4 and by our assumption, a contradiction.

Thus after all $\bar{w}=0.$ So
\begin{eqnarray*}
1 &=&\lim_{\theta \rightarrow \infty }[e^{\bar{\mu}\theta }\Phi (\bar{W}%
(\theta )-\sigma )+e^{\sigma \bar{W}(\theta )}\Phi (0)] \\
&=&\lim_{\theta \rightarrow \infty }\frac{e^{\frac{1}{2}\bar{\sigma}%
^{2}\theta }e^{-\frac{1}{2}\bar{W}^{2}-\frac{1}{2}\bar{\sigma}^{2}\theta
+\sigma \bar{W}(\theta )}}{(\sigma -\bar{W}(\theta ))\sqrt{2\pi }}+e^{\sigma
\bar{W}(\theta )}\Phi (0)] \\
&=&\lim_{\theta \rightarrow \infty }e^{\sigma \bar{W}(\theta )}[\frac{1}{2}-%
\frac{1}{\sigma \sqrt{2\pi }}]=\frac{1}{2}\lim_{\theta \rightarrow \infty
}e^{\sigma \bar{W}(\theta )}.\qquad \square
\end{eqnarray*}
\bigskip

\noindent \textbf{Proof of Conclusion 1:} For any $\varepsilon ,$ put
\begin{equation*}
W_{\varepsilon }(\theta ):=\frac{\bar{\mu}-\frac{1}{2}\bar{\sigma}%
^{2}+\varepsilon }{\bar{\sigma}}\sqrt{\theta }:\qquad e^{-\bar{\mu}\theta
+\sigma W_{\varepsilon }(\theta )+\frac{1}{2}\sigma ^{2}}=e^{\varepsilon
\sqrt{\theta }}.
\end{equation*}%
For $\varepsilon >0$ and large enough $\theta ,$
\begin{eqnarray*}
e^{-\bar{\mu}\theta -\sigma \bar{W}(\theta )+\frac{1}{2}\sigma ^{2}} &=&\Phi
(-\bar{W}(\theta ))+\varphi (-\bar{W}(\theta ))/H(\sigma -\bar{W}(\theta ))
\\
&<&e^{\varepsilon \sqrt{\theta }}=e^{-\bar{\mu}\theta +\sigma W_{\varepsilon
}(\theta )+\frac{1}{2}\sigma ^{2}}:\qquad -\bar{W}(\theta )<W_{\varepsilon
}(\theta ).
\end{eqnarray*}%
On the other hand, for $\varepsilon <0$ and large enough $\theta $
\begin{equation*}
e^{-\bar{\mu}\theta -\sigma \bar{W}(\theta )+\frac{1}{2}\sigma ^{2}}>\Phi (-%
\bar{W}(\theta ))>e^{\varepsilon \sqrt{\theta }}=e^{-\bar{\mu}\theta +\sigma
W_{\varepsilon }(\theta )+\frac{1}{2}\sigma ^{2}}:\qquad W_{\varepsilon
}(\theta )<-\bar{W}(\theta ).\qquad \square
\end{equation*}
\bigskip

\noindent \textbf{Proof of Lemma 6. }Here, for large $\theta $,
\begin{equation*}
\bar{W}(\theta )=-(\bar{\mu}-\frac{1}{2}\bar{\sigma}^{2})\sqrt{\theta }+o(%
\sqrt{\theta }).
\end{equation*}%
So
\begin{equation*}
\sigma +\bar{W}(\theta )=\frac{\frac{3}{2}\bar{\sigma}^{2}-\bar{\mu}}{\bar{%
\sigma}}\sqrt{\theta }+o(\sqrt{\theta }),\qquad \sigma -\bar{W}(\theta )=%
\frac{\frac{1}{2}\bar{\sigma}^{2}+\bar{\mu}}{\bar{\sigma}}\sqrt{\theta }+o(%
\sqrt{\theta }).
\end{equation*}%
Furthermore, rewriting the normal censor equation,%
\begin{equation*}
e^{-\bar{\mu}\theta +\frac{1}{2}\sigma ^{2}-\sigma \bar{W}(\theta )}=e^{-%
\frac{1}{2}\bar{W}(\theta )^{2}}\Phi (-\sigma +\bar{W}(\theta ))/\varphi (-%
\bar{W}(\theta )+\sigma )+\Phi (-\bar{W}(\theta )).
\end{equation*}%
So by Lemma 4, and since $\bar{W}(\theta )\rightarrow -\infty ,$
\begin{equation*}
\lim_{\theta \rightarrow \infty }e^{-\bar{\mu}\theta +\frac{1}{2}\sigma
^{2}-\sigma \bar{W}(\theta )}=1.
\end{equation*}%
In fact, we have $e^{-\bar{\mu}\theta +\frac{1}{2}\sigma ^{2}-\sigma \bar{W}%
(\theta )}=1+o(1/\sqrt{\theta }).$ $\square $
\bigskip

\noindent \textbf{Proof of Conclusion 2. }As $2\bar{\mu}>\bar{\sigma}^{2},$
note that
\begin{equation*}
(\bar{\sigma}^{2}-\bar{\mu})-\frac{1}{2}(\frac{3}{2}\bar{\sigma}-\frac{\bar{%
\mu}}{\bar{\sigma}})^{2}=(\bar{\sigma}^{2}-\bar{\mu})-\frac{9}{8}\bar{\sigma}%
^{2}-\frac{1}{2}\frac{\bar{\mu}^{2}}{\bar{\sigma}^{2}}+\frac{3}{2}\bar{\mu}=-%
\frac{1}{8}\bar{\sigma}^{2}+\frac{1}{2}\bar{\mu}-\frac{1}{2}\frac{\bar{\mu}%
^{2}}{\bar{\sigma}^{2}}<0;
\end{equation*}%
indeed,%
\begin{equation*}
\bar{\mu}^{2}-\bar{\mu}\bar{\sigma}^{2}+\frac{1}{4}\bar{\sigma}^{4}=(\bar{\mu%
}-\frac{1}{2}\bar{\sigma}^{2})^{2}>0.
\end{equation*}%
From \S 6,%
\begin{eqnarray*}
g(\theta ) &=&e^{(\bar{\sigma}^{2}-\bar{\mu})\theta }\Phi (\sigma +\bar{W}%
(\theta ))+\Phi (-\bar{W}(\theta ))+o(1/\sqrt{\theta }) \\
&=&e^{(\bar{\sigma}^{2}-\bar{\mu})\theta }\Phi (\frac{\frac{3}{2}\bar{\sigma}%
^{2}-\bar{\mu}}{\bar{\sigma}}\sqrt{\theta })+\Phi (\frac{\bar{\mu}-\frac{1}{2%
}\bar{\sigma}^{2}}{\bar{\sigma}}\sqrt{\theta })+o(1/\sqrt{\theta }).
\end{eqnarray*}
Applying the asymptotic expansion \cite{AbrS}
\begin{equation*}
\Phi (x)\sim 1-\frac{e^{-x^{2}/2}}{x\sqrt{2\pi }}\qquad (\text{as }%
x\rightarrow +\infty ),
\end{equation*}%
yields%
\begin{equation*}
e^{(\bar{\sigma}^{2}-\bar{\mu})\theta }\Phi (\frac{\frac{3}{2}\bar{\sigma}%
^{2}-\bar{\mu}}{\bar{\sigma}}\sqrt{\theta })=e^{(\bar{\sigma}^{2}-\bar{\mu}%
)\theta }+o(1/\sqrt{\theta }).\qquad \square
\end{equation*}

\noindent Mathematics Department, University of Leicester, University Road, Leicester LE1 7RH;
\newline
Mathematics Department, London School of Economics, Houghton Street, London
WC2A 2AE; 
\newline 
A.J.Ostaszewski@lse.ac.uk

\end{document}